\documentclass[11pt]{amsproc}
\usepackage{amsmath,amsthm,amssymb,amsfonts,hyperref,enumerate}
\usepackage[usenames,dvipsnames,svgnames,table]{xcolor} 
\hypersetup{
 colorlinks=true,
 citecolor=Blue,
 linkcolor=Red,
 urlcolor=Violet}
\title[Representations of automorphisms on holomorphic differentials]{Explicit Galois representations of automorphisms on holomorphic differentials in characteristic $p$}
\author{Kenneth A. Ward}
\address{Kenneth A. Ward \\ NYU-ECNU Institute of Mathematical Sciences, NYU Shanghai, People's Republic of China 200122. Ph: +86(21) 2059-5124. Fax: +1(212) 995-4195. \\ E-mail: kennethward@nyu.edu}

\newtheorem{thm}{Theorem}
\newtheorem{lem}{Lemma}
\newtheorem{cor}{Corollary}
\allowdisplaybreaks
\begin{document}
\maketitle

\begin{abstract}
We determine the representation of the group of automorphisms for cyclotomic function fields in characteristic $p > 0$ induced by the natural action on the space of holomorphic differentials via construction of an explicit basis of differentials. This includes those cases which present wild ramification and automorphism groups with non-cyclic $p$-part, which have remained elusive. We also obtain information on the gap sequences of ramified primes. Finally, we extend these results to rank one Drinfel'd modules. \\ 

\end{abstract}

\section{Introduction}
Let $K$ be an algebraic function field with algebraically closed constant field $k$ and genus $g_K \geq 2$. The space $\Omega_K$ of holomorphic differentials of $K$ is a vector space over $k$, and represents the finite group of automorphisms $G$ of $K/k$, as introduced by Hurwitz \cite{Hur}. If char $k=0$, the representation of $G$ induced by $\Omega_K$ was completely determined by Chevalley and Weil \cite{ChWeHe}. In char $k > 0$, this problem remains open. Various methods have been introduced to solve this problem with certain assumptions on ramification or group structure, but in general the structure of $\Omega_K$ as a representation space for $G$ has been unknown. In this section, we review the previous work on this problem and summarize our approach to an explicit solution in positive characteristic. 

The construction of Chevalley and Weil may be completed algebraically by evaluating the different of the extension $K/K^G$. This extends to all cases in char $k = p > 0$ where $(|G|,p) = 1$ \cite{VaMa2}. In positive characteristic, Tamagawa had resolved the unramified cyclic case \cite{Tam}, using Hasse and Witt's generation of unramified extensions of degree $p$ to show that $\Omega_K \cong I_G \oplus (k[G])^{g_{K^G}-1}$, where $I_G$ denotes the identity representation of $G$. This remains valid for unramified extensions in characteristic zero \cite{GlGrHe}. Tamagawa conjectured that this was likely true for all non-cyclic extensions in positive characteristic, which is not correct. Using the Frattini subgroup, Hasse-Witt theory, and relative projective $k[G]$-modules, Valentini was able to show for unramified covers that Tamagawa's decomposition is equivalent to the existence of a cyclic $p$-Sylow subgroup $P$ so that its normalizer $N$ may be written as a direct product of $P$ with $N/P$ \cite{Val}. It is key for this argument that the cover be unramified, as differentials are then not altered in passing from fixed fields of $p$-Sylow subgroups to $K$.

Via Kani and Nakajima, in the presence of ramification, we have for tame covers the Grothendieck equivalence \begin{equation}\label{eq1} [\Omega_K] \cong [ k \oplus (k[G])^{g_{K^G}-1} \oplus \tilde{R}_G^* ], \end{equation} where exactness is derived via the sequence $$0 \rightarrow \Omega_K \rightarrow \Omega_K(S) \xrightarrow{Res} kS \xrightarrow{\Sigma} k \rightarrow 0,$$ with $\Omega_K(S)$ denoting the differentials of $K$ logarithmic along $S$, and $\tilde{R}_G^*$ the contragredient of a unique module $\tilde{R}_G$ satisfying $\tilde{R}_G^{\oplus |G|} = R_G$, where $R_G$ denotes the ramification module of $G$ \cite{Kan,Nak1,Nak2}. The condition on ramification is essential for this argument: By Serre duality and the Riemann-Hurwitz formula, $\Omega_K(S)$ is a projective $k[G]$-module if, and only if, ramification is tame. As in Valentini's result, one may show that we may avoid passing to Grothendieck groups in \eqref{eq1} if, and only if, $G$ is $p$-nilpotent and possesses cyclic $p$-Sylow subgroups. Other results are possible when $p \mid |G|$ due to calculations of Boseck \cite{Bos,Gar1}, which assume that $G$ is cyclic or elementary abelian, or that ramified primes are totally ramified \cite{Gar2,Kako,Madden,RzViMa1,RzViMa2,VaMa2}. We will omit a complete discussion of these, but these are all motivated by, and often directly following, Boseck's construction. If $G$ possesses non-cyclic $p$-part, difficulties with characterisation of indecomposable $k[G]$-modules obstruct a group-theoretic formulation, as Higman showed that in this case $G$ has indecomposable representations of arbitrarily large degree \cite{Hig}. 

We are interested in an approach which circumvents all of these obstructions. As Rzedowski-Calder\'{o}n et al. have mentioned, the proofs in many previously known cases are similar, despite that the arithmetic objects under study are distinct \cite{RzViMa2}. It is particularly necessary to consider situations where the automorphism group presents wild ramification and non-cyclic $p$-part, and where ramification is not total at the ramified primes. A natural and broad class of such objects is given to us by cyclotomic function fields as formulated by Carlitz \cite{Car}, where ramification is understood via explicit class field theory \cite{Hay1,Hay2}. Cyclotomic function fields are of special interest to us due to their analogy to the cyclotomic extensions of $\mathbb{Q}$. For the Frobenius $\phi_q(u):= u^q$ ($q = p^r$), $\mu_T:=Tu$, and $M=\sum_{k=0}^n a_k T^k \in \mathbb{F}_q[T]$, the Carlitz polynomial is defined as \begin{equation}\label{eq2} u_q^M := \sum_{k=0}^n a_k (\phi_q + \mu_T)^k(u),\end{equation} and the Carlitz-Hayes module as $\Lambda_{q,M} = \{u \; |\; u_q^M = 0\}$, within some fixed choice of algebraic closure of $\mathbb{F}_q(T)$. We denote the cyclotomic function field $\mathbb{F}_q(T)(\Lambda_{q,M})$ by $K_{q,M}$ and its automorphism group by $G_{q,M}$. Links to cyclotomic fields in the classical case are numerous: For example, an isomorphism \begin{equation}\label{eq3} G_{q,M} \cong (\mathbb{F}_q[T]/(M))^*\end{equation} is given by the $\mathbb{F}_q[T]$-module action \eqref{eq2}, the composite of $\overline{\mathbb{F}}_q$ and all cyclotomic extensions of $\mathbb{F}_q(T)$ and $\mathbb{F}_q(1/T)$ (viewing the variable as $X = 1/T$ in the latter case to capture wild ramification at infinity) is equal to the maximal abelian extension of $\mathbb{F}_q(T)$, and the $\mathbb{F}_q$-module isomorphism $\phi_P: \Lambda_{q,M} \rightarrow \Lambda_{q,M}$ defined as $\phi_P(\lambda) = \lambda^P$ for irreducible $P \nmid M$ ($P \in \mathbb{F}_q[T]$) is given by the Artin map \cite{Hay1}. If the factorisation of $M$ contains at least one square, then $G_{q,M}$ generally exhibits both non-cyclic $p$-part and wild ramification. The reader is referred elsewhere for a complete description of automorphism group structure and ramification in cyclotomic function fields \cite{Ros,Vil}. As the field $K_{q,M}$ is geometric over $\mathbb{F}_q(T)$, we may pass to an algebraic closure $\overline{\mathbb{F}}_q$, preserving \eqref{eq3} as well as ramification in $K_{q,M}/\mathbb{F}_q(T)$. The constant field is commonly assumed to be algebraically closed in order to guarantee that gap sequences and Weierstrass points, which behave differently in positive characteristic than in characteristic zero, are well-defined \cite{GaVi,Schmid,Schmidt2}. Here, we shall work over $\mathbb{F}_q$, assuming only that the degree of each ramified prime is equal to one, which is equivalent to the splitting of $M$ in $\mathbb{F}_q$. We shall refer to these as the \emph{split} cyclotomic function fields. We point out that all of our calculations remain valid over $\overline{\mathbb{F}}_q$. 

Our first main result identifies an explicit canonical basis over $\mathbb{F}_q$ for the holomorphic differentials, which we will denote by $\Omega_{q,M}$, of $K_{q,M}$. This basis is expressed in terms of generators of Carlitz-Hayes modules, and is localized at an arbitrary choice of ramified prime of $\mathbb{F}_q[T]$. Our second main result describes the representation of $G_{q,M}$ induced by the explicit basis of $\Omega_{q,M}$, which is decomposed locally at each irreducible factor $P$ of $M$ via $P$-adic forms of elements of $G_{q,M}$ according to \eqref{eq3}. As a corollary, we identify gap sequences of totally ramified primes. Finally, we generalize our arguments to Drinfel'd modules of rank one and explain why our results are the best possible.

An outline of the subsequent sections of this paper is as follows. In Section 2, we find the $\mathbb{F}_q$-basis of $\Omega_{q,M}$. We describe the corresponding representation of $G_{q,M}$ in Section 3. Section 4 addresses questions on gap sequences. In Section 5, we provide the extension to Drinfel'd modules and discuss the case where ramified primes are of degree greater than one. We conclude with some acknowledgements.

\section{The basis of $\Omega_{q,M}/\mathbb{F}_q$}
Methods of Boseck and others rely upon modulo $p$ reductions. Using calculations of the different, we construct bases of $\Omega_{q,M}/\mathbb{F}_q$ via $P$-adic expansions for irreducible polynomials $P \in \mathbb{F}_q[T]$. The ramified primes of $K_{q,M}/\mathbb{F}_q(T)$ correspond to the factors of $M$ and the prime at infinity, and the primes of $K_{q,M}$ above infinity are of degree one for every cyclotomic function field \cite[Corollary 1.2]{GaRo}. Henceforth, we assume without loss that $M$ and its factors are monic. We begin with $M$ square-free. \\

\noindent {\sc Case 1: $M$ is square-free}. Let $M = \prod_{i=1}^r P_i$ for distinct, irreducible $P_i \in \mathbb{F}_q[T]$ ($i=1,\ldots,r$). If $r=1$, this is trivial, as then $g_{q,M}=0$ and $K_{q,M}$ possesses no holomorphic differentials other than zero. If $r = 2$, then $M = P_1 P_2$ for distinct, irreducible $P_1, P_2 \in \mathbb{F}_q[T]$. Let $\mathfrak{p}_1$ and $\mathfrak{p}_2$ denote the primes of $\mathbb{F}_q(T)$ associated with $P_1$ and $P_2$, respectively, and let $\mathfrak{p}_\infty$ again denote the prime of $\mathbb{F}_q(T)$ at infinity. Each prime $\mathfrak{p}_1$, $\mathfrak{p}_2$, and $\mathfrak{p}_\infty$ is ramified in $K_{q,M}$ with ramification index equal to $q-1$, and the primes of $K_{q,P_1}$ above $\mathfrak{p}_1$ and $\mathfrak{p}_\infty$ are unramified in $K_{q,M}$ (see \cite[Proposition 2.2]{Hay1}). Thus the different $\mathfrak{D}_{K_{q,M}/K_{q,P_1}}$ of $K_{q,M}$ over $K_{q,P_1}$ is equal to \begin{equation}\label{eq4} \mathfrak{D}_{K_{q,M}/K_{q,P_1}} = \prod_{\mathfrak{P}|\mathfrak{p}_2} \mathfrak{P}^{q-2}.\end{equation} By \eqref{eq4} and the Riemann-Hurwitz formula, the genus $g_{q,M}$ of $K_{q,M}$ is given by \begin{align}\label{eq5}\notag g_{q,M} &= 1 + [K_{q,M}:K_{q,P_1}](g_{q,P_1}-1) + \frac{1}{2} d_{K_{q,M}}(\mathfrak{D}_{K_{q,M}/K_{q,P_1}}) \\&= 1 - (q-1) + \frac{1}{2}(q-1)(q-2) \\&\notag = \frac{(q-3)(q-2)}{2}, \end{align} where $d_{K_{q,M}}(\cdot)$ denotes the degree function on ideals of $K_{q,M}$. The different $\mathfrak{D}_{q,M}$ of $K_{q,M}/\mathbb{F}_q(T)$ is equal to $$\mathfrak{D}_{q,M} = \prod_{\mathfrak{P}|\mathfrak{p}_1} \mathfrak{P}^{q-2} \prod_{\mathfrak{P}|\mathfrak{p}_2} \mathfrak{P}^{q-2} \prod_{\mathfrak{P}|\mathfrak{p}_\infty} \mathfrak{P}^{q-2}.$$ This implies that the differential $dT$ has divisor in $K_{q,M}$ equal to \begin{align} \label{eq6} \notag (dT)_{K_{q,M}} & = \mathfrak{D}_{K_{q,M}} \text{Con}_{K_{q,M}/\mathbb{F}_q(T)}(\mathfrak{p}_\infty^{-2}) \\& = \left(  \prod_{\mathfrak{P}|\mathfrak{p}_1} \mathfrak{P}^{q-2} \prod_{\mathfrak{P}|\mathfrak{p}_2} \mathfrak{P}^{q-2} \prod_{\mathfrak{P}|\mathfrak{p}_\infty} \mathfrak{P}^{q-2} \right) \prod_{\mathfrak{P}|\mathfrak{p}_\infty} \mathfrak{P}^{-2(q-1)} \\& \notag =  \prod_{\mathfrak{P}|\mathfrak{p}_1} \mathfrak{P}^{q-2} \prod_{\mathfrak{P}|\mathfrak{p}_2} \mathfrak{P}^{q-2} \prod_{\mathfrak{P}|\mathfrak{p}_\infty} \mathfrak{P}^{-q}. \end{align} For each $i=1,2$, the generator $\lambda_i$ of $\Lambda_{q,P_i}$ over $\mathbb{F}_q[T]$ has divisor in $K_{q,M}$ equal to \begin{equation} \label{eq7} (\lambda_i)_{K_{q,M}} = \frac{\prod_{\mathfrak{P}|\mathfrak{p}_i} \mathfrak{P}}{\prod_{\mathfrak{P}|\mathfrak{p}_\infty} \mathfrak{P}}.\end{equation} It follows from \eqref{eq6} and \eqref{eq7} that $$(\lambda_1^{-\mu_1} \lambda_2^{-\mu_2} dT)_{K_{q,M}} = \prod_{\mathfrak{P}|\mathfrak{p}_1} \mathfrak{P}^{q-2-\mu_1} \prod_{\mathfrak{P}|\mathfrak{p}_2} \mathfrak{P}^{q-2-\mu_2} \prod_{\mathfrak{P}|\mathfrak{p}_\infty} \mathfrak{P}^{\mu_1 + \mu_2 - q}.$$ Thus $\lambda_1^{-\mu_1} \lambda_2^{-\mu_2} dT \in \Omega_{q,M}$ if $\mu_1,\mu_2 \leq q-2$ and $\mu_1 + \mu_2 \geq q$. Therefore the number of such differentials is equal to $\sum_{k=2}^{q-2} (k-1) = g_{q,M}$. We may define the set of possible values of $\mu_1$ and $\mu_2$ as $$\Gamma_2(P_1) = \{(\mu_1,\mu_2) \in \mathbb{Z}^2 \; | \; \mu_1 \leq q-2,\;0 \leq \mu_2 \leq q-2,\; \mu_1 + \mu_2 \geq q\}.$$ As $[K_{q,M}:K_{q,P_1}] = q-1$, the differentials $\lambda_1^{-\mu_1} \lambda_2^{-\mu_2} dT$ with $(\mu_1,\mu_2) \in \Gamma_2(P_1)$ are linearly independent over $\mathbb{F}_q$. As $|\Gamma_2(P_1)| = g_{q,M}$, these differentials must form a basis of $\Omega_{q,M}/\mathbb{F}_q$. 

Inductively, let $r \geq 3$. We suppose that $M = \prod_{i=1}^{r} P_i$, for distinct irreducible $P_i \in \mathbb{F}_q[T]$. Let $\lambda$ denote a generator of $\Lambda_{q,M}$ over $\mathbb{F}_q[T]$. Via the prime decomposition $\Lambda_{q,M} = \oplus_{i=1}^r \Lambda_{q,P_i}$, the generator $\lambda$ may be uniquely written as $\lambda=\sum_{i=1}^r \lambda_i$, where $\lambda_i$ denotes a generator of $\Lambda_{q,P_i}$ over $\mathbb{F}_q[T]$, for each $i=1,\ldots,r$. We let $F = \prod_{i=1}^{r-1} P_i$. By assumption, the elements \begin{equation*} \prod_{i=1}^{r-1} \lambda_i^{-\mu_i} dT \;\;\;\;\;\;\;\;\;\; \left(\mu_1 \leq q-2,\;0 \leq \mu_2,\ldots,\mu_{r-1}\leq q-2,\sum_{i=1}^{r-1} \mu_i \geq q\right)\end{equation*} form a basis of $\Omega_{q,F}/\mathbb{F}_q$. If $\sum_{i=1}^{r-1} \mu_i \geq q$, then by previous arguments, a differential of the form $\prod_{i=1}^r \lambda_i^{-\mu_i} dT$ lies in $\Omega_{q,M}$ if $0 \leq \mu_r \leq q-2$. Otherwise, if $2 \leq \sum_{i=1}^{r-1} \mu_i = k < q$, the number of possible values for $\mu_r$ is $k-1$, and as $\mu_r \leq q-2$, there exist no such holomorphic differentials of $K_{q,M}$ with $\sum_{i=1}^{r-1} \mu_i < 2$. We define the set of differentials \begin{align*} \mathfrak{T}_{q,M}(\lambda;P_1) = \Bigg\{  \prod_{i=1}^{r} \lambda_i^{-\mu_i} dT \;\bigg|&\; \mu_1 \leq q-2, \; \\& 0 \leq \mu_2,\ldots,\mu_{r} \leq q-2,\; \sum_{i=1}^{r} \mu_i \geq q \Bigg\}.\end{align*} By construction, $\mathfrak{T}_{q,M}(\lambda;P_1) \subset \Omega_{q,M}$. Also, we define \begin{align}\label{eq25} \Gamma_{r-1,k}(P_1) = \Bigg\{& (\mu_1,\ldots,\mu_{r-1}) \in \mathbb{Z}^{r-1} \;\bigg| \\&\notag \mu_1 \leq q-2, \; 0 \leq \mu_2,\ldots,\mu_{r-1} \leq q-2, \; \sum_{i=1}^{r-1} \mu_i = k\Bigg\}.\end{align} We find that the number of elements in the set $\mathfrak{T}_{q,M}(\lambda;P_1)$ is equal to \begin{align}\label{eq8} |\mathfrak{T}_{q,M}(\lambda;P_1)| & = (q-1)\sum_{k=q}^\infty |\Gamma_{r-1,k}(P_1)| + \sum_{k=2}^{q-1} (k-1) |\Gamma_{r-1,k}(P_1)| \\& \notag = (q-1)g_{q,F} + \sum_{k=2}^{q-1} (k-1)|\Gamma_{r-1,k}(P_1)|. \end{align} Let $k \leq q-1$. By induction on $r$, it follows that 
   \begin{equation*}|\Gamma_{r-1,k}(P_1)| = \begin{cases}
       (q-1)^{r-2} & \text{if }2 \leq k < q-1\\
       (q-1)^{r-2}-1 & \text{if }k=q-1.
     \end{cases}
\end{equation*} The prime $\mathfrak{p}_{r}$ of $\mathbb{F}_q(T)$ associated with $P_{r}$ is unramified in $K_{q,F}/\mathbb{F}_q(T)$, each prime of $K_{q,F}$ dividing $\mathfrak{p}_{r}$ is totally and tamely ramified in $K_{q,M}/K_{q,F}$, and all other primes are unramified in $K_{q,M}/K_{q,F}$. Thus $d_{K_{q,M}}(\mathfrak{D}_{K_{q,M}/K_{q,F}}) = (q-1)^{r-1} (q-2)$. By the Riemann-Hurwitz formula, we obtain that \begin{align} \label{eq9} \notag (q-1)g_{q,F}& + \sum_{k=2}^{q-1}(k-1)|\Gamma_{r-1,k}(P_1)| \\& \notag = (q-1)g_{q,F} + \sum_{k=1}^{q-2} k|\Gamma_{r-1,k+1}(P_1)| \\& = (q-1)g_{q,F} + \left( \sum_{k=1}^{q-2} k(q-1)^{r-2} \right) - (q-2) \\& \notag = 1 + (q-1)(g_{q,F}-1) + \frac{1}{2}(q-1)^{r-1}(q-2) \\& \notag = 1 + [K_{q,M}:K_{q,F}](g_{q,F}-1) + \frac{1}{2}d_{K_{q,M}}(\mathfrak{D}_{K_{q,M}/K_{q,F}}) \\& \notag = g_{q,M}. \end{align} The elements of $\mathfrak{T}_{q,M}(\lambda;P_1)$ are linearly independent over $\mathbb{F}_q$ by construction. By \eqref{eq8} and \eqref{eq9}, $|\mathfrak{T}_{q,M}(\lambda;P_1)| = g_{q,M}$. Therefore $\mathfrak{T}_{q,M}(\lambda;P_1)$ forms a basis of $\Omega_{q,M}/\mathbb{F}_q$. For each $i=1,\ldots,r$, we now define \begin{align}\label{eq35} \mathfrak{B}_{q,M}(\lambda;P_i) = \Bigg\{& \prod_{j=1}^r \lambda_j^{-\mu_j} P_i^{\mu_0} dT \; \bigg| \;\mu_0 \geq 0, \\&\notag 0 \leq \mu_1,\ldots,\mu_r \leq q-2, \; \sum_{j=1}^r \mu_j - (q-1)\mu_0 \geq q\Bigg\}.\end{align} The reason for the notation $\mathfrak{B}_{q,M}(\lambda;P_i)$ will be made clear later. As $\lambda_i^{q-1} = -P_i$ by definition of $\lambda_i$ and the $\mathbb{F}_q[T]$-module action \eqref{eq2}, this is consistent with the definition of $\mathfrak{T}_{q,M}(\lambda;P_1)$, up to multiplication of elements by $\pm 1$. By induction, we have shown that if $\deg (M) \geq 2$, then $\mathfrak{B}_{q,M}(\lambda;P_i)$ forms a basis of $\Omega_{q,M}/\mathbb{F}_q$. We have thus proven the following lemma.

\begin{lem} Let $M = \prod_{i=1}^r P_i \in \mathbb{F}_q[T]$ $(r \geq 2)$ be square-free with each $P_i \in \mathbb{F}_q[T]$ linear. Then for each $i=1,\ldots,r$, the set $\mathfrak{B}_{q,M}(\lambda;P_i)$ forms a basis of $\Omega_{q,M}/\mathbb{F}_q$. 
\end{lem}

We shall now proceed to treat wild ramification in $K_{q,M}/\mathbb{F}_q(T)$, which occurs if $M$ is not square-free. We accomplish this by first addressing the case where $M$ is a power of an irreducible polynomial. \\

\noindent {\sc Case 2: $M = P^n$ ($n \geq 2$)}. Let $\mathfrak{P}$ denote the prime of $K_{q,M}$ lying above the prime $\mathfrak{p}$ of $\mathbb{F}_q(T)$ associated with $P$. Via explicit class field theory, we recognize the different $\mathfrak{D}_{q,M}$ of $K_{q,M}/\mathbb{F}_q(T)$ as \begin{equation}\label{eq30}\mathfrak{D}_{q,M} = \mathfrak{P}^s \prod_{\mathfrak{A}|\mathfrak{p}_\infty} \mathfrak{A}^{q-2}, \end{equation} where $s = nq^n - (n+1)q^{n-1}$ \cite[Theorem 4.1]{Hay1}. Let $\lambda_n$ be a generator of $\Lambda_{q,M}$ over $\mathbb{F}_q[T]$. For each $i=1,\ldots,n$, the map $\theta(\lambda_n) = \lambda_n^{P^{n-i}}$ induces an isomorphism of $\mathbb{F}_q$-modules $\Lambda_{q,P^n} / \Lambda_{q,P^{n-i}} \cong \Lambda_{q,P^i}$, so that the element $\lambda_i = \lambda_n^{P^{n-i}}$ generates $\Lambda_{q,P^i}$ over $\mathbb{F}_q[T]$. The prime $\mathfrak{p}$ is totally ramified in $K_{q,M}$ and for the unique prime $\mathfrak{P}_i$ of $K_{q,P^i}$ above $\mathfrak{p}$, $v_{\mathfrak{P}_i}(\lambda_i)=1$. Therefore the valuation of $(\prod_{i=2}^n \lambda_i^{\mu_i})\lambda_1^{-\mu_1}dT$ at $\mathfrak{P}$ is equal to \begin{equation} \label{eq22} nq^n - (n+1)q^{n-1} - q^{n-1} \mu_1 + \sum_{i=2}^n q^{n-i} \mu_i.\end{equation} As the prime $\mathfrak{p}_\infty$ is unramified in $K_{q,M}/K_{q,P}$, $v_{\mathfrak{P}_\infty}(\lambda_1) = -1$ at the prime $\mathfrak{P}_\infty$ of $K_{q,P}$ above $\mathfrak{p}_\infty$, and $(dT)_{K_{q,M}} = \mathfrak{D}_{q,M} \text{Con}_{K_{q,M}/\mathbb{F}_q(T)} (\mathfrak{p}_\infty^{-2}) = \mathfrak{P}^s \prod_{\mathfrak{A}|\mathfrak{p}_\infty} \mathfrak{A}^{-q}$, the valuation of $\lambda_1^{-\mu_1} dT$ at any prime $\mathfrak{A}$ of $K_{q,M}$ lying above $\mathfrak{p}_\infty$ is equal to $\mu_1-q$. As in the proof of the Brumer-Stark conjecture for global function fields, for a certain prime $\mathfrak{A}_\infty$ of $K_{q,M}$ above $\mathfrak{p}_\infty$, we may select the generator $\lambda_n$ to satisfy \begin{equation} \label{eq23} v_{\mathfrak{A}_\infty}(\lambda_n^A) = (q-1)(n-\deg(A) - 1) - 1 \end{equation} for each $A \in \mathbb{F}_q[T]\backslash \{0\}$ with $\deg(A) < n$ \cite[Theorem 12.14]{Ros}. Let $\mathbb{F}_q(T)_\infty$ denote the completion of $\mathbb{F}_q(T)$ at $\mathfrak{p}_\infty$, and let $\iota: K_{q,M} \rightarrow \overline{\mathbb{F}_q(T)}_\infty$ be an embedding which corresponds to $\mathfrak{A}_\infty$. For each $A \in \mathbb{F}_q[T]$ with $\deg (A) < \deg (M)$, we find by choice of $\lambda_n$ that $v_\infty((\iota(\lambda_n))^A) = n - \deg (A) - 1 - (q-1)^{-1}.$ For each $i=1,\ldots,n$ and $A \in \mathbb{F}_q[T]$ with $\deg (A) < i$, we thus obtain \begin{align*} v_\infty(\iota(\lambda_i^A))& =v_\infty(\iota((\lambda_n^{P^{n-i}})^A))\\& = v_\infty((\iota(\lambda_n))^{P^{n-i} A}) \\& = n - \deg (P^{n-i} A) - 1 - (q-1)^{-1} \\& = i - \deg(A) - 1 - (q-1)^{-1}. \end{align*} Also, for each $i=1,\ldots,n$ and $A \in (\mathbb{F}_q(T)/(P^i))^*$, the map $\sigma_A(\lambda_i) = \lambda_i^A$ induces an automorphism of $K_{q,P^i}/\mathbb{F}_q(T)$. We conclude for such $A$ that \begin{align} \label{eq10} \notag v_{\sigma_A^{-1}(\mathfrak{A}_\infty|_{K_{q,P^i}})} (\lambda_i) &= v_{\mathfrak{A}_\infty|_{K_{q,P^i}}}(\lambda_i^A) \\& = (q-1)(i-\deg(A) - 1) - 1 \\& \notag \geq -1. \end{align} Via the isomorphism $G_{q,P^i} \cong (\mathbb{F}_q(T)/(P^i))^*$, the automorphisms of $K_{q,P^i}$ are given by $\sigma_A$ for $A \in (\mathbb{F}_q(T)/(P^i))^*$, and these act transitively on the primes of $K_{q,P^i}$ above $\mathfrak{p}_\infty$. Furthermore, every prime of $K_{q,P^i}$ above $\mathfrak{p}_\infty$ is unramified in $K_{q,M}$. It follows for any prime $\mathfrak{A}$ of $K_{q,M}$ above $\mathfrak{p}_\infty$ that $v_\mathfrak{A} (\lambda_i) \geq -1$, and for such $\mathfrak{A}$ that $v_\mathfrak{A}(\left( \prod_{i=2}^n \lambda_i^{\mu_i}\right) \lambda_1^{-\mu_1} dT) \geq \mu_1 - \sum_{i=2}^n \mu_i - q$. Therefore $(\prod_{i=2}^n \lambda_i^{\mu_i}) \lambda_1^{-\mu_1} dT \in \Omega_{q,M}$, provided that $\mu_1,\ldots,\mu_n$ are chosen so that \begin{equation}\label{eq11} nq^n - (n+1)q^{n-1} - q^{n-1}\mu_1 + \sum_{i=2}^n q^{n-i}\mu_i \geq 0,\;\;\;\;\;\; \mu_1-\sum_{i=2}^n \mu_i - q \geq 0.\end{equation} Let us also suppose that $\mu_1 \leq nq - (n+1)$ and $\mu_i \geq 0$ for each $i=2,\ldots,n$. If $n=2$, such differentials are those which satisfy $q \leq \mu_1 \leq 2q-3$ and $0 \leq \mu_2 \leq \mu_1-q$. The genus of $K_{q,M}$ is given by $g_{q,M} = 1 + \frac{1}{2}(q^2 - 3q)$, and the number of these differentials is equal to $\sum_{k=q}^{2q-3} (k-q+1) = g_{q,M}$. As $0 \leq \mu_2 \leq q-3$ for any value of $\mu_1$, these differentials are linearly independent over $\mathbb{F}_q$, and hence form a basis of $\Omega_{q,M}/\mathbb{F}_q$. If $n \geq 3$, we require that $0 \leq \mu_i \leq q-1$ for each $i=2,\ldots,n$ to guarantee linear independence over $\mathbb{F}_q$. For all integers $n \geq 2$, we define \begin{align*} \notag \mathfrak{W}_{q,M}(\lambda_n;P) = \Bigg\{&\left(\prod_{i=2}^n \lambda_i^{\mu_i}\right) \lambda_1^{-\mu_1} dT\;\bigg|\;\\& nq^n - (n+1)q^{n-1} - q^{n-1}\mu_1 + \sum_{i=2}^n q^{n-i}\mu_i \geq 0, \\&\notag \mu_1 - \sum_{i=2}^n \mu_i - q \geq 0,\; \mu_1 \leq nq-(n+1),\\&\notag 0 \leq \mu_i \leq q-1,\; i=2,\ldots,n\Bigg\}.\end{align*} We have already shown that the set $\mathfrak{W}_{q,M}(\lambda_n;P)$ forms a basis of $\Omega_{q,M}/\mathbb{F}_q$ in the case $n=2$. We shall demonstrate that this is also valid for all $n \geq 3$. We may view the power series \begin{align}\label{eq31} \left(\sum_{i=0}^{q-1} x^i \right)^{n-1} & = \left[\sum_{l=0}^{n-1} (-1)^l {n-1 \choose l} x^{ql}\right]\left[\sum_{t=0}^\infty {n-2+t \choose t} x^t \right] = \sum_{m=0}^\infty a_m x^m \end{align} as a generating function for $\mathfrak{W}_{q,M}(\lambda_n;P)$ (see for instance \cite{Mur}). Via \eqref{eq34} and \eqref{eq31}, we find that if $n \geq 3$, then the number of elements in the set $\mathfrak{W}_{q,M}(\lambda_n;P)$ is equal to \begin{align} \label{eq12} \sum_{k=0}^{(n-1)q-(n+1)} &\sum_{m=0}^k \sum_{l=0}^{\left\lfloor \frac{m}{q}\right\rfloor} (-1)^l {n-1 \choose l} {n-2 + m - ql \choose n - 2}  \\&\notag= \sum_{l=0}^{\left\lfloor n-1-\frac{n+1}{q}\right\rfloor} \;\; \sum_{k=0}^{(n-1-l)q - (n+1)} (-1)^l {n-1 \choose l} \sum_{m=n-2}^{k+n-2} {m \choose n-2}, \end{align} where $\lfloor \cdot \rfloor$ denotes the floor function. The inner sum in the last expression of the equation \eqref{eq12} is equal to the $(n-2)$nd coefficient of the polynomial \begin{equation} \label{eq13} \sum_{m=n-2}^{k+n-2} (1+X)^m = X^{-1}[(1+X)^{k+n-1} - (1+X)^{n-2}].\end{equation} It follows from \eqref{eq12} and \eqref{eq13} that $|\mathfrak{W}_{q,M}(\lambda_n;P)|$ is equal to the $n$th coefficient of the polynomial \begin{equation}\label{eq27} \sum_{l=0}^{\left\lfloor n-1-\frac{n+1}{q}\right\rfloor} (-1)^l {n-1 \choose l}(1+X)^{(n-1-l)q-1}.\end{equation} The $n$th coefficient of \eqref{eq27} is the same as that of $(1+X)^{-1}[((1+X)^q-1)^{n-1}+(-1)^n]$, which may be computed directly as \begin{align*} q^{n-1}\Bigg[\frac{1}{2}(q-1)&(n-1)-1\Bigg] + 1  \\& = 1 - [K_{q,M}:\mathbb{F}_q(T)] + \frac{1}{2}[(nq^n - (n+1)q^{n-1}) + q^{n-1}(q-2)] \\& = 1 + [K_{q,M}:\mathbb{F}_q(T)](g_{\mathbb{F}_q(T)} - 1) + \frac{1}{2} d_{K_{q,M}}(\mathfrak{D}_{q,M})\\& = g_{q,M}. \end{align*} Thus $|\mathfrak{W}_{q,M}(\lambda_n;P)| = g_{q,M}$. Also, the elements of $\mathfrak{W}_{q,M}(\lambda_n;P)$ are linearly independent over $\mathbb{F}_q$ by construction. Therefore $\mathfrak{W}_{q,M}(\lambda_n;P)$ forms a basis of $\Omega_{q,M}/\mathbb{F}_q$. We now define for all integers $n \geq 2$ the set \begin{align} \label{eq34}\mathfrak{B}_{q,M}(\lambda_n;P) = \Bigg\{&\left(\prod_{i=2}^n \lambda_i^{\mu_i}\right) \lambda_1^{-\mu_1}P^{\mu_0} dT\;\bigg|\;\\&\notag nq^n - (n+1)q^{n-1} - q^{n-1}\mu_1 + \sum_{i=2}^n q^{n-i}\mu_i \geq 0, \\&\notag \mu_1 - \sum_{i=2}^n \mu_i - q \geq 0,\;\mu_0 \geq 0,\\&\notag (n-1)(q-1) \leq \mu_1 \leq nq-(n+1), \\&\notag 0 \leq \mu_i \leq q-1,\; i=2,\ldots,n\Bigg\}.\end{align} The reason for the notation $\mathfrak{B}_{q,M}(\lambda_n;P)$ will once again be made clear later. As $\lambda_1^{q-1} = -P$ by definition of $\lambda_1$ and the $\mathbb{F}_q[T]$-module action \eqref{eq2}, this is consistent with the definition of $\mathfrak{W}_{q,M}(\lambda_n;P)$, up to multiplication of elements by $\pm 1$. We have thus proven the following lemma. 
 
\begin{lem} Let $M = P^n \in \mathbb{F}_q[T]$ $(n \geq 2)$ with $P \in \mathbb{F}_q[T]$ linear. Then the set $\mathfrak{B}_{q,M}(\lambda_n;P)$ forms a basis of $\Omega_{q,M}/\mathbb{F}_q$. \end{lem}

We now suppose that $M = \prod_{i=1}^r P_i^{n_i}$ for distinct, irreducible $P_i \in \mathbb{F}_q[T]$. For each $i=1,\ldots,r$, we let $\mathfrak{P}_i$ denote any prime of $K_{q,M}$ above the prime $\mathfrak{p}_i$ of $\mathbb{F}_q(T)$, where $\mathfrak{p}_i$ is associated with $P_i$. As $\mathfrak{p}_i$ is unramified in $K_{q,M}$ outside of $K_{q,P_i^{n_i}}$, the bound on valuation for required for holomorphicity of a differential of $K_{q,M}$ at $\mathfrak{P}_i$ is precisely the same as the analogous bound at the prime of $K_{q,P_i^{n_i}}$ above $\mathfrak{p}_i$. On the other hand, at a prime of $K_{q,M}$ above $\mathfrak{p}_\infty$, the corresponding bound on valuation is an amalgam of the Kummer and wildly ramified cases, as we will see momentarily. For each $k=1,\ldots,n_j$ and $j=1,\ldots,r$, we let $\lambda_{j,k}$ denote a generator over $\mathbb{F}_q[T]$ of $\Lambda_{q,P_j^k}$, which is chosen to satisfy $\lambda_{j,n_j}^{P^{n_j - k}} = \lambda_{j,k}$. Via the prime decomposition $\Lambda_{q,M} = \oplus_{i=1}^r \Lambda_{q,P^{n_i}}$, the element $\lambda = \sum_{i=1}^r \lambda_{i,n_i}$ is a generator of the $\mathbb{F}_q[T]$-module $\Lambda_{q,M}$, and this decomposition of $\lambda$ is unique. For each $i=1,\ldots,r$, we thus define \begin{align}\label{eq33}\mathfrak{B}_{q,M} &(\lambda; P_i) = \Bigg\{ \prod_{j=1}^r \left[\left( \prod_{k=2}^{n_j} \lambda_{j,k}^{\mu_{j,k}} \right) \lambda_{j,1}^{-\mu_{j,1}} \right] P_i^{\mu_0} \; dT \; \bigg| \\&\notag\;\;\;\;\;\;  n_i q^{n_i} - (n_i+1)q^{n_i - 1} + q^{n_i-1}(q-1)\mu_0 - q^{n_i - 1} \mu_{i,1} + \sum_{k=2}^{n_i} q^{n_i - k} \mu_{i,k} \geq 0, \\&\notag\;\;\;\;\;\; n_j q^{n_j} - (n_j+1)q^{n_j - 1} - q^{n_j - 1} \mu_{j,1} + \sum_{k=2}^{n_j} q^{n_j - k} \mu_{j,k} \geq 0, \; j \neq i, \\&\notag\;\;\;\;\;\; \sum_{j=1}^r \mu_{j,1} - (q-1)\mu_0 - \sum_{j=1}^r \left( \sum_{k=2}^{n_j} \mu_{j,k} \right) - q \geq 0, \; \mu_0 \geq 0, \\&\notag\;\;\;\;\;\; (n_j - 1)(q-1) \leq \mu_{j,1} \leq n_j q - (n_j + 1),\\&\notag\;\;\;\;\;\; 0 \leq \mu_{j,k} \leq q-1, \; k \geq 2, \; j = 1,\ldots,r \Bigg\}.\end{align} This is consistent with our previous definitions \eqref{eq35} and \eqref{eq34} of $\mathfrak{B}_{q,M}(\cdot)$ in Lemmas 1 and 2. By construction, the set $\mathfrak{B}_{q,M}(\lambda;P_i)$ consists of $\mathbb{F}_q$-linearly independent holomorphic differentials. 

In order to show that the cardinality of $\mathfrak{B}_{q,M}(\lambda;P_i)$ is equal to the genus $g_{q,M}$, we may employ induction and the methods of the proofs of Lemmas 1 and 2. We assume without loss that $i=1$, and also that $n_1 > 1$ or $r > 1$ (or else there are no holomorphic differentials other than zero). If $r = 1$, it is then immediate from Lemma 2 that $|\mathfrak{B}_{q,M}(\lambda;P_1)| = |\mathfrak{B}_{q,P_1^{n_1}}(\lambda_{1,n_1};P_1)| = g_{q,M}$. If $r>1$, then in analogy to $\Gamma_{r-1,k}(P_1)$ in the proof of Lemma 1 (see \eqref{eq25}), we may define $\Phi_{r-1,k}(P_1)$ for each integer $k$ as the set of $\mu = (\mu_0,\mu_{1,1},...,\mu_{r-1,n_{r-1}}) \in \mathbb{Z} \oplus (\oplus_{j=1}^{r-1} \mathbb{Z}^{n_j})$ which satisfy the equality \begin{equation*} \sum_{j=1}^{r-1} \mu_{j,1} - (q-1)\mu_0 - \sum_{j=1}^{r-1} \left(\sum_{l=2}^{n_j} \mu_{j,l} \right) = k,\end{equation*} in addition to all of the inequalities for $j = 1, \ldots, r-1$ which appear in the definition of $\mathfrak{B}_{q,M}(\lambda;P_1)$. Each element of $\omega \in \mathfrak{B}_{q,M}(\lambda;P_1)$ corresponds to some $\mu \in \Phi_{r-1,k}(P_1)$, for some integer $k$. We shall denote this correspondence by $\omega \sim \mu$, and we let $\nu_{r,k}(\mu) = |\{\omega \in \mathfrak{B}_{q,M}(\lambda;P_1)\;|\; \omega \sim \mu\}|$. The quantity $\nu_{r,k}(\mu)$ depends only on $r$ and $k$, and not on the choice of $\mu$, as interaction between valuations above distinct primes of $\mathbb{F}_q(T)$ in the definition of $\mathfrak{B}_{q,M}(\lambda;P_1)$ occurs only at infinity. We thus let $\nu_{r,k} = \nu_{r,k}(\mu)$ for any choice of $\mu \in \Phi_{r-1,k}(P_1)$. If $n_r = 1$, then as in the proof of Lemma 1 (see the argument preceding \eqref{eq25}), $\nu_{r,k}=\min\{\max\{0,k-1\},q-1\}$. If $n_r \geq 2$, then similarly to \eqref{eq31}, we may view the Laurent series \begin{align} \label{eq29}&\;\;\;\;\;x^{-[n_r q-(n_r+1)]}\left(\sum_{i=0}^{q-2} x^i\right) \left(\sum_{j=0}^{q-1} x^j\right)^{n_r-1} \\&\notag=x^{-[n_r q-(n_r+1)]}(1-x^{q-1})\left[\sum_{l=0}^{n_r-1} (-1)^l {n_r-1 \choose l} x^{ql}\right]\left[\sum_{t=0}^\infty {n_r-1+t \choose t} x^t \right]\\&\notag = \sum_{m \in \mathbb{Z}} a_{r,m} x^m\end{align} as a generating function corresponding to the possible values of $\mu_{r,1},\ldots,\mu_{r,n_r}$. By definition \eqref{eq33}, $\nu_{r,k} = \sum_{m \leq k-q} a_{r,m}$. The coefficients $a_{r,m}$ may be calculated as in \eqref{eq31} in the proof of Lemma 2. For the inductive step, if $r=2$, then $|\Phi_{r-1,k}(P_1)|=|\Phi_{1,k}(P_1)|$ may be found directly from \eqref{eq34} and \eqref{eq31} via $$|\Phi_{1,k}(P_1)| = \sum_{m=k}^{n_1 q - (n_1+1)} a_{m-k}.$$ If $r > 2$, then we may calculate $|\Phi_{r-1,k}(P_1)|$ using \eqref{eq29} via $$|\Phi_{r-1,k}(P_1)|=\sum_l a_{r-1,l}|\Phi_{r-2,k+l}(P_1)|.$$ As in the proof of Lemma 1 (see \eqref{eq8}, \eqref{eq9}), we may then invoke the Riemann-Hurwitz formula for the extension $K_{q,M}/\mathbb{F}_q(T)$ \cite[Theorem 12.7.2]{Vil} to obtain \begin{equation*}|\mathfrak{B}_{q,M}(\lambda;P_1)| = \sum_k \nu_{r,k} |\Phi_{r-1,k}(P_1)| = \sum_k |\Phi_{r-1,k}(P_1)| \sum_{m \leq k - q} a_{r,m} = g_{q,M}.\end{equation*} The arguments follow the proofs of Lemmas 1 and 2, to which we refer the reader for further details. We have thus obtained the following theorem. 

\begin{thm} Suppose that $M \in \mathbb{F}_q[T]$ is of degree at least two and splits in $\mathbb{F}_q$. Then the set $\mathfrak{B}_{q,M}(\lambda;P_i)$ forms a basis of $\Omega_{q,M}/\mathbb{F}_q$ (and the condition on the degree of $M$ is necessary for $\Omega_{q,M} \neq \{0\}$). \end{thm}

We reiterate that Theorem 1 is equally valid over $\overline{\mathbb{F}}_q$. We shall refer to the set $\mathfrak{B}_{q,M}(\lambda;P_i)$ as a \emph{canonical basis of $\Omega_{q,M}$ at $P_i$}. We note that a canonical basis of $\Omega_{q,M}$ at $P_i$ may indeed be described only in terms of $\lambda$, the differential $dT$, and the $\mathbb{F}_q[T]$-module action \eqref{eq2}.

\section{The representation of $G_{q,M}$}

The representation of $G_{q,M}$ afforded by its action on $\Omega_{q,M}$ may be described via $P$-adic decompositions for each irreducible factor $P$ of $M$. If $M = P^n$, the action of $G_{q,M}$ on $\Omega_{q,M}$ may be completely described in terms of the canonical basis $\mathfrak{B}_{q,M}(\lambda;P)$ of $\Omega_{q,M}$. The action of $G_{q,M}$ on this basis is an explicitly defined linear transformation, as we will see below. If $M$ is square-free, as in the Kummer case, elements of $G_{q,M}$ act as roots of unity on differentials, according to their decompositions as elements of $(\mathbb{F}_q[T]/(M))^*$ via \eqref{eq3}. Excluding these two cases, the action of $G_{q,M}$ on a canonical basis does not admit a natural expression in terms of the $\mathbb{F}_q[T]$-module action. We shall find that this may be resolved. 

We again let $M = \prod_{i=1}^r P_i^{n_i}$. With the notation of Section 2, we define the finite set \begin{align} \label{eq24} \notag\mathfrak{V}_{q,M} (\lambda) \;=\; &\Bigg\{ \prod_{i=1}^r \left[\left( \prod_{k=2}^{n_i} \lambda_{i,k}^{\mu_{i,k}} \right) \lambda_{i,1}^{-\mu_{i,1}} \right] dT \; \bigg|\; \\&\; \sum_{i=1}^r \mu_{i,1} - \sum_{i=1}^r \left( \sum_{k=2}^{n_i} \mu_{i,k} \right) - q \geq 0, \\\notag&\; n_i q^{n_i} - (n_i + 1) q^{n_i - 1} - q^{n_i - 1} \mu_{i,1} + \sum_{k=2}^{n_i} q^{n_i - k} \mu_{i,k} \geq 0, \\&\notag\; \mu_{i,1} \leq n_i q - (n_i+1),\\&\notag\; 0 \leq \mu_{i,k} \leq q-1,\; k \geq 2, \; i=1,\ldots,r \Bigg\}. \end{align} We have already shown that $\bigcup_{i=1}^r \mathfrak{B}_{q,M}(\lambda;P_i) \subset \mathfrak{V}_{q,M}(\lambda) \subset \Omega_{q,M}$. In particular, $\mathfrak{V}_{q,M}(\lambda)$ generates $\Omega_{q,M}$ over $\mathbb{F}_q$. \\

\noindent {\sc Step 1: The $\mathbb{F}_q[T]$-module action}. Via the isomorphism \eqref{eq3}, for each $A \in (\mathbb{F}_q[T]/(M))^*$ and $i=1,\ldots,r$, we may associate $\sigma_A \in G_{q,M}$ with the $P_i$-adic decomposition of $A$, \begin{equation} \label{eq16} A = \sum_{l=0}^{\deg(A)} \alpha_{i,l} P_i^l \;\;\;\;\;\; (\alpha_{i,l} \in \mathbb{F}_q, \; l=0,\ldots,\deg(A)).\end{equation} The action of $\sigma_A$ on the local component $(\prod_{k=2}^{n_i} \lambda_{i,k}^{\mu_{i,k}})\lambda_{i,1}^{-\mu_{i,1}}$ in the notation of \eqref{eq24} for an element of $\mathfrak{V}_{q,M}(\lambda)$ is given according to \eqref{eq16} and the definition of the $\mathbb{F}_q[T]$-module action \eqref{eq2} by the $P_i$-adic form \begin{equation} \label{eq17} \sigma_A(\lambda_{i,k}) = \lambda_{i,k}^A = \lambda_{i,k}^{\sum_{l=0}^{\deg(A)} \alpha_{i,l} P_i^l} = \sum_{l=0}^{\deg(A)} \alpha_l \lambda_{i,k}^{P_i^l} = \sum_{l=0}^{k-1} \alpha_l \lambda_{i,k-l}. \end{equation} At the prime $\mathfrak{P}_i$ of $K_{q,M}$ above $\mathfrak{p}_i$, we have $v_{\mathfrak{P}_i}(\lambda_{i,k-l}) = q^{n_i-(k-l)} = q^l q^{n_i-k} = q^l v_{\mathfrak{P}_i}(\lambda_{i,k})$. At a prime $\mathfrak{A}$ of $K_{q,M}$ above $\mathfrak{p}_\infty$, the valuations of $\lambda_{i,k}$ and $\lambda_{i,k-l}$ are given by \eqref{eq10}. It follows from this and \eqref{eq17} that, in the notation of \eqref{eq24}, the action of $\sigma_A$ on an element of $\mathfrak{V}_{q,M}(\lambda)$ is equal to an $\mathbb{F}_q$-linear combination of elements of the form \begin{equation}\label{eq18} \omega = \prod_{i=1}^r \left[\left( \prod_{k=2}^{n_i} \lambda_{i,k}^{a_{i,k}} \right) \lambda_{i,1}^{-a_{i,1}} \right] dT \in \Omega_{q,M},\end{equation} each of which satisfies the bounds \begin{align}\label{eq19}\notag & \mathfrak{P}_i\; (i=1,\ldots,r):\; n_i q^{n_i} - (n_i+1)q^{n_i-1} - q^{n_i-1} a_{i,1} + \sum_{k=2}^{n_i} q^{n_i - k} a_{i,k} \geq 0,  \\& \quad\quad\quad\quad\quad\quad\quad\quad\; a_{i,1} \leq n_i q - (n_i+1), \\& \notag \mathfrak{A}|\mathfrak{p}_\infty: \; \sum_{i=1}^r a_{i,1} - \sum_{i=1}^r \left( \sum_{k=2}^{n_i} a_{i,k} \right) - q \geq 0. \end{align} However, the differential \eqref{eq18} does not necessarily lie in $\mathfrak{V}_{q,M}(\lambda)$, as at least one of the inequalities $0 \leq a_{i,k} \leq q-1$ may not hold for some $i=1,\ldots,r$ and $k \geq 2$. We see in the following step that this may be resolved.\\

\noindent {\sc Step 2: Reduction to $\mathfrak{V}_{q,M}(\lambda)$}. By definition of the $\mathbb{F}_q[T]$-module action on $\Lambda_{q,M}$, we find for $k \geq 2$ that \begin{equation} \label{eq20} \lambda_{i,k}^q = \lambda_{i,k}^{P_i} - P_i \lambda_{i,k} = \lambda_{i,k-1} - P_i \lambda_{i,k} = \lambda_{i,k-1} + \lambda_{i,1}^{q-1} \lambda_{i,k}.\end{equation} Let $m=2,\ldots,n_i$ (provided that $n_i \geq 2$) for which $a_{i,m} \geq q$. Using \eqref{eq20}, we may thus write $\omega$ of \eqref{eq18} as \begin{equation} \label{eq21} \omega = \prod_{i=1}^r \Bigg[\Bigg( \prod_{k=2}^{n_i} \lambda_{i,k}^{b_{i,k}} \Bigg) \lambda_{i,1}^{-b_{i,1}} \Bigg] dT + \prod_{i=1}^r \left[\left( \prod_{k=2}^{n_i} \lambda_{i,k}^{c_{i,k}} \right) \lambda_{i,1}^{-c_{i,1}} \right] dT,\end{equation} where $a_{j,k} = b_{j,k} = c_{j,k}$ for all $k$ and $j \neq i$, and $a_{i,k} = b_{i,k} = c_{i,k}$ for all $k \neq 1, m-1,m$. Also, for $m>2$, $$\begin{array}{ll} b_{i,1} = a_{i,1}, & c_{i,1} = a_{i,1} - (q-1),\\ b_{i,m-1} = a_{i,m-1} + 1, & c_{i,m-1} = a_{i,m-1}, \\ b_{i,m} = a_{i,m}-q, & c_{i,m} = a_{i,m} - (q-1), \end{array}$$ and similarly if $m=2$. We therefore obtain the bounds at $\mathfrak{P}_i$ ($i=1,\ldots,r$) and $\mathfrak{A}|\mathfrak{p}_\infty$ for the first term in the decomposition of $\omega$ in \eqref{eq21} if $m >2$ (and similarly if $m=2$): \begin{flalign*} &\mathfrak{P}_i:\;  n_i q^{n_i} - (n_i+1) q^{n_i-1} - q^{n_i-1}b_{i,1} + \sum_{k=2}^{n_i} q^{n_i - k} b_{i,k} & \\& \quad\quad=\; q^{n_i - m}(a_{i,m} - q) + q^{n_i-(m-1)}(a_{i,m-1}+1) + n_i q^{n_i} - (n_i+1)q^{n_i-1} \\& \quad\quad\quad\quad\quad\quad\quad\quad\quad\quad\quad\quad\quad\quad\quad\quad\quad\quad\quad\;\; -q^{n_i-1} a_{i,1} + \sum_{\substack{k \geq 2 \\ k \neq m-1, m}} q^{n_i-k}a_{i,k}  \\& \quad\quad=\; n_i q^{n_i} - (n_i+1) q^{n_i - 1} - q^{n_i-1} a_{i,1} + \sum_{k=2}^{n_i} q^{n_i - k} a_{i,k} \geq 0, \\& \quad\quad\; b_{i,1} = a_{i,1} \leq n_i q - (n_i+1), \\& \mathfrak{A}|\mathfrak{p}_\infty: \; \sum_{i=1}^r b_{i,1} - \sum_{i=1}^r \left( \sum_{k=2}^{n_i} b_{i,k} \right) - q &\\&\quad\quad\quad\; = \; -(a_{i,m} - q) - (a_{i,m-1} + 1) + a_{i,1} - \sum_{\substack{k \geq 2 \\ k \neq m-1,m}} a_{i,k} + \sum_{j \neq i} a_{j,1} \\& \quad\quad\quad\quad\quad\quad\quad\quad\quad\quad\quad\quad\quad\quad\quad\quad\quad\quad\quad\quad\quad\quad\quad - \sum_{j \neq i} \left( \sum_{k=2}^{n_j} a_{j,k} \right) - q \\&\quad\quad\quad\; = \; q-1 + \sum_{i=1}^r a_{i,1} - \sum_{i=1}^r \left(\sum_{k=2}^{n_i} a_{i,k} \right) - q \geq q-1.\end{flalign*} We also find analogous bounds for the second term in the decomposition of $\omega$ in \eqref{eq21}: \begin{flalign*}&\mathfrak{P}_i:\; n_i q^{n_i} - (n_i+1)q^{n_i-1} - q^{n_i - 1} c_{i,1} + \sum_{k=2}^{n_i} q^{n_i - k}c_{i,k} \geq \; (q^{n_i-1} - q^{n_i - m})(q-1), & \\& \quad\quad\; c_{i,1} = a_{i,1} -(q-1)\leq n_i q - (n_i+1)-(q-1), \\ & \mathfrak{A}|\mathfrak{p}_\infty : \;  \sum_{i=1}^r c_{i,1} - \sum_{i=1}^r \left(\sum_{k=2}^{n_i} c_{i,k} \right) - q \geq 0.& \end{flalign*} Furthermore, the valuations of all differentials in \eqref{eq21} at all other primes of $K_{q,M}$ are equal to zero. Therefore the first and second terms on the right-hand side of \eqref{eq21} both also lie in $\Omega_{q,M}$. 

Inductive application of Step 2, via descent within $\omega$ from $\lambda_{i,n_i}$ ($i=1,\ldots,r$) in the first and second terms on the right-hand side of \eqref{eq21}, permits elements of the form $\omega$ as in \eqref{eq18} to be written in the manner prescribed as $\mathbb{Z}$-sums, hence over $\mathbb{F}_p$, of elements of $\mathfrak{V}_{q,M}(\lambda)$. By Steps 1 and 2, the action of $\sigma_A$ on an element of $\mathfrak{V}_{q,M}(\lambda)$ thus admits explicit description as an $\mathbb{F}_q$-linear combination of elements of $\mathfrak{V}_{q,M}(\lambda)$. This yields the following theorem.

\begin{thm} Suppose that $M \in \mathbb{F}_q[T]$ is of degree at least two and splits in $\mathbb{F}_q$. Then the action of $G_{q,M}$ on $\Omega_{q,M}$ is determined by its $\mathbb{F}_q[T]$-module action on $\mathfrak{V}_{q,M}(\lambda)$, and may be explicitly described in terms of $P_i$-adic forms of elements of $G_{q,M}$ via the isomorphism \eqref{eq3}. \end{thm}

If $M=P^n$, we find that $\mathfrak{V}_{q,M}(\lambda)=\mathfrak{B}_{q,M}(\lambda;P)$, which gives the following corollary to Theorem 2.

\begin{cor} Suppose that $M = P^n \in \mathbb{F}_q[T]$ is of degree at least two and splits in $\mathbb{F}_q$, with $P$ irreducible. Then the action of $G_{q,M}$ on $\Omega_{q,M}$ is determined explicitly by its $\mathbb{F}_q[T]$-module action on the canonical basis $\mathfrak{B}_{q,M}(\lambda;P)$ according to $P$-adic forms of elements of $G_{q,M}$ via the isomorphism \eqref{eq3}. \end{cor}

Henceforth, we shall refer to $\mathfrak{V}_{q,M}$ as a set of \emph{canonical generators of $\Omega_{q,M}$}.

\section{Gap sequences of $K_{q,M}$}

Identification of spaces of holomorphic differentials can be useful for determining gap sequences. For example, Garcia showed that for most Artin-Schreier extensions of $k(x)$ (where $k$ is algebraically closed), the ramified primes are Weierstrass points, by determining the gap sequences for the ramified primes via examination of orders of holomorphic differentials in Boseck's basis \cite{Gar1}. Other examples include the Hermite bases for holomorphic differentials of the curve $y^q + y = x^{q+1}$: $$\{(x-a)^s P^r dx \;|\;r,s \geq 0,\; r + s \leq q-2\},$$ where $P = (y-b) - a^q (x-a)$ and $(a,b)$ is a finite point on the curve. As shown by Garcia and Viana, the Weierstrass points for this curve correspond to the rational points over $\mathbb{F}_{q^2}$, and are identified as such via determination of the order sequence at each finite point on the curve using the associated Hermite basis \cite[Theorem 2]{GaVi}. 

As all of our calculations over $\mathbb{F}_q$ carry through to $\overline{\mathbb{F}}_q$, we may consider gap sequences for the split cyclotomic function fields. With the notation of Section 2, we consider the case of $M = P^n$, and we let $\Gamma(P^n)$ denote the set of all tuples $(\mu_1,\ldots,\mu_n)$ appearing as respective powers of $\lambda_1,\ldots,\lambda_n$ in an element of the canonical basis $\mathfrak{B}_{q,M}(\lambda_n;P)$. The formula \eqref{eq22} makes explicit the valuations of such differentials at the prime $\mathfrak{P}$ of $K_{q,M}$ above $\mathfrak{p}$. As the valuation at $\mathfrak{P}$ of each element of $\mathfrak{B}_{q,M}(\lambda_n;P)$ is distinct in this case, we obtain the following corollary to Theorem 1.

\begin{cor} Suppose that $M = P^n \in \mathbb{F}_q[T]$ is of degree at least two and splits in $\mathbb{F}_q$, with $P$ irreducible. Then the order sequence of $K_{q,M}$ at $\mathfrak{P}$ is given by \begin{equation*} \left\{nq^n - (n+1)q^{n-1} - q^{n-1}\mu_1 + \sum_{i=2}^n q^{n-i}\mu_i\;\bigg|\; (\mu_1,\ldots,\mu_n) \in \Gamma(P^n)\right\}. \end{equation*} \end{cor}

If the automorphism group is cyclic and of prime power order, it is possible to determine that the totally ramified primes are Weierstrass points \cite{VaMa1}. This type of argument employs Witt vector decompositions to establish rapid growth of the degree of the different at the totally ramified primes. For cyclotomic function fields, the automorphism group is generally neither cyclic nor of prime power order, and the degree of the different at ramified primes does not exhibit sufficient growth to apply such methods \cite{Hay1}. 

\section{Extensions}

These results extend to Drinfel'd modules of rank one in a natural way. Let $K$ denote any algebraic function field with constant field $\mathbb{F}_q$, $A$ the ring of elements of $K$ with pole only at a prime $\mathfrak{P}_\infty$ of $K$, $\mathfrak{m}$ a non-zero proper ideal of $A$, and $\rho$ a rank one Drinfel'd $A$-module over $\mathbb{F}_q$. We also let $H^+_A$ denote a narrow class field, i.e., a normalizing field for rank one Drinfel'd $A$-modules over $(K,\mathfrak{P}_\infty,\text{sgn})$ for some choice of $\text{sgn}$ map, and $K(\mathfrak{m}):=H_A^+(\rho[\mathfrak{m}])$ \cite{Dri,Hay3}. The extension $K(\mathfrak{m})/H_A^+$ is abelian and exhibits all of the arithmetic properties of the cyclotomic function fields, including ramification. Furthermore, the $A$-module $\rho[\mathfrak{m}]$ is cyclic with $|(A/\mathfrak{m})^*|$ generators. We thus obtain the following corollary to our main results in Sections 3 and 4. For simplicity, we assume that $\mathfrak{m} = \prod_{i=1}^r \mathfrak{p}_i^{n_i}$ with each $\mathfrak{p}_i$ distinct and irreducible, as well as all other notations analogous to those of previous sections.

\begin{cor} Suppose that $\mathfrak{m}$ is of degree at least two and splits in $\mathbb{F}_q$. Then one may define for the space of holomorphic differentials $\Omega_{q,\mathfrak{m}}$ a canonical basis $\mathfrak{B}_{q,\mathfrak{m}}(\lambda;\mathfrak{p}_i)$ at $\mathfrak{p}_i$ and generating set $\mathfrak{V}_{q,\mathfrak{m}}(\lambda)$. The action of $G_{q,\mathfrak{m}}:=\text{\emph{Aut}}(K(\mathfrak{m})/H_A^+)$ on $\Omega_{q,\mathfrak{m}}$ is given by the $\mathfrak{p}_i$-adic decompositions of elements of $G_{q,\mathfrak{m}}$ via the isomorphism $G_{q,\mathfrak{m}} \cong (A/\mathfrak{m})^*$.
\end{cor}

In the same manner that an algebraically closed constant field guarantees well-defined gap sequences, the case where $M$ does not split in $\mathbb{F}_q$ behaves quite differently from the split cyclotomic function fields. For example, for a non-linear irreducible $P \in \mathbb{F}_q[T]$, the zero divisor of a generator $\lambda$ of $\Lambda_{q,P}$ over $\mathbb{F}_q[T]$ contains primes above infinity. Explicitly, with $d = \deg(P)$, one may show that the zero divisor of $\lambda$ in $K_{q,P}$ contains precisely $(q^{d-1}-1)/(q-1)$ such primes, where the valuations of $\lambda$ at these primes are determined according to \eqref{eq23}. In fact, valuations are not sufficiently controlled in this case to allow for the same constructions of canonical bases completed in Section 2, and our results are thus the best possible.

\section{Acknowledgements}

The author would like to thank A. Ghosh and B.-H. Im for their helpful comments and questions, as well as NYU Shanghai and the NYU-ECNU Institute of Mathematical Sciences for their support.

\bibliographystyle{plain}
\raggedright
\bibliography{references}
\end{document}